\documentclass[a4]{amsart}
\numberwithin{equation}{section}
\newtheorem{theorem}{Theorem}[section]
\newtheorem{lemma}[theorem]{Lemma}
\newtheorem{corollary}[theorem]{Corollary}
\DeclareMathOperator{\modu}{Mod}
\DeclareMathOperator{\Hom}{Hom}
\DeclareMathOperator{\Ext}{Ext}
\DeclareMathOperator{\cok}{cok}

\DeclareMathOperator{\ddim}{\underline{dim}}
\DeclareMathOperator{\Rep}{Rep}
\DeclareMathOperator{\add}{add}
\DeclareMathOperator{\proj}{proj}
\DeclareMathOperator{\Supp}{Supp}
\DeclareMathOperator{\Sl}{Sl}
\DeclareMathOperator{\Tr}{Tr}
\DeclareMathOperator{\Gl}{Gl}
\keywords{Invariant theory, quivers}
\subjclass{Primary 13A50}
\author{Aidan Schof\hspace{.070757em}ield}
\address{Department of Mathematics\\ University of Bristol
Senate House \\
Tyndall Avenue \\
Bristol BS8 1TH
}
\email{Aidan.Schofield@bristol.ac.uk}
\author{Michel Van den Bergh}
 \address{Departement WNI\\Limburgs
  Universitair Centrum\\ Universitaire Campus\\ Building D\\ 3590
  Diepenbeek\\ Belgium} 
\email{vdbergh@luc.ac.be,
  http://www.luc.ac.be/Research/Algebra/} 
\thanks{The second author is
  a senior researcher at the FWO}
\date{\today}
\title{Semi--invariants of quivers for arbitrary dimension vectors}
\begin{document}
\begin{abstract}
The representations of dimension vector $\alpha$ of the quiver $Q$ can be
parametrised by a vector space $R(Q,\alpha)$ on which an algebraic
group $\Gl(\alpha)$ acts so that the set of orbits is bijective with
the set of isomorphism classes of representations of the quiver. We
describe the semi--invariant polynomial functions on this vector space
in terms of the category of representations. More precisely, we
associate to a suitable map between projective representations a
semi--invariant polynomial function that describes when this map is
inverted on the representation and we show that these semi--invariant
polynomial functions form a spanning set of all semi--invariant
polynomial functions in characteristic $0$. If the quiver has no
oriented cycles, we may replace consideration of inverting maps between
projective representations by consideration of representations that
are left perpendicular to some representation of dimension vector
$\alpha$. These left perpendicular representations are just the
cokernels of the maps between projective representations that we
consider.
\end{abstract} 
\maketitle

\section{Notation and Introduction}
\label{s0}
In the sequel $k$ will be an algebraically closed field. For our
main result, Theorem \ref{t2} below, $k$ will have characteristic zero.

Let $Q$ be a quiver with finite vertex set $V$, finite arrow set $A$
and two functions $i,t:A\rightarrow V$ where for an arrow $a$ we shall
usually write $ia$ in place of $i(a)$, the initial vertex, and
$ta=t(a)$, the terminal vertex. A representation, $R$, of the quiver
$Q$ associates a $k$-vector space $R(v)$ to each vertex $v$ of the quiver
and a linear map $R(a):R(ia)\rightarrow R(ta)$ to each arrow $a$. A
homomorphism $\phi$ of representations from $R$ to $S$ is given by a
collection of linear maps for each vertex $\phi(v):R(v)\rightarrow
S(v)$ such that for each arrow $a$, $R(a)\phi(ta)=\phi(ia)S(v)$. The
category of representations of the quiver $\Rep(Q)$ is an abelian
category as is the full subcategory of finite dimensional
representations. We shall usually be interested in finite dimensional
representations in which case each representation has a dimension
vector $\ddim R$, which is a function from the set of vertices $V$ to
the natural numbers $\mathbb{N}$ defined by $(\ddim R)(v) = \dim
R(v)$. Now let $\alpha$ be a dimension vector for $Q$. The
representations of dimension vector $\alpha$ are parametrised by the
vector space
\begin{equation*}
R(Q,\alpha) = \times_{a \in A}{^{\alpha(ia)}k^{\alpha(ta)}}
\end{equation*}
where ${^{m}k^{n}}$ is the vector space of $m$ by $n$ matrices over
$k$ (and ${^{m}k}$ and $k^{n}$ are shorthand for ${^{m}k^{1}}$ and
${^{1}k^{n}}$ respectively). Given a point $p \in R(Q,\alpha)$, we
denote the corresponding representation by $R_{p}$. The isomorphism
classes of representations of dimension vector $\alpha$ are in $1$ to
$1$ correspondence with the orbits of the algebraic group
\begin{equation*}
\Gl(\alpha) = \times_{v \in V} \Gl_{\alpha(v)}(k)
\end{equation*}
We consider the action of $\Gl(\alpha)$ 
on the co--ordinate ring $S(Q,\alpha)$ of $R(Q,\alpha)$; 
$f \in S(Q,\alpha)$ is said to be semi--invariant of weight $\psi$ 
where $\psi$ is a character of $\Gl(\alpha)$ if 
$g(f) = \psi(g)f, \; \forall g \in \Gl(\alpha)$.

The invariants and semi--invariants for this action are of
importance for the description of the moduli spaces of representations
of the quiver for fixed dimension vector. In characteristic $0$ we may
apply Weyl's theory of invariants for $\Sl_{n}(k)$ to give an explicit
description of all such semi--invariants. We shall find a set of
semi--invariants that span all semi--invariants as a vector space in
characteristic $0$. It seems likely that the result we obtain here
should hold in arbitrary characteristic and that this would follow
from Donkin \cite{donkin,donkin1}. However, we restrict ourselves in this paper to
characteristic $0$. In \cite{schofield}, the first author described all the
semi--invariants when $\Gl(\alpha)$ has an open orbit on $R(Q,\alpha)$
in terms of certain polynomial functions naturally associated to the
representation theory of the quiver. We shall begin by recalling the
definition of these semi--invariants and some related theory.

Given the quiver $Q$, let $\add(Q)$ be the additive $k$--category
generated by $Q$. To describe this more precisely, we define a path
of length $n>0$ from the vertex $v$ to the vertex $w$ to be a monomial
in the arrows $p=a_{1}\dots a_{n}$ such that $ta_{i}=ia_{i+1}$ for
$0<i<n$ and $ia_{1}=v$, $ta_{n}=w$. We define $ip=v$ and $tp=w$. For
each vertex $v$ there is also the trivial path of length $0$, $e_{v}$,
from $v$ to $v$.

For each vertex $v$, we have an object $O(v)$ in $\add(Q)$,
and $\Hom(O(v),O(w)) = P(v,w)$ where $P(v,w)$ is the vector space with
basis the paths from $v$ to $w$ including the trivial path if $v = w$;
finally,
\begin{equation*}
\Hom(\bigoplus_{v}O(v)^{a(v)}, \bigoplus_{v}O(v)^{b(v)})
\end{equation*}
is defined in the usual way for an additive category, where
composition arises via matrix multiplication.

Any representation $R$ of $Q$ extends uniquely to a covariant 
functor from $\add(Q)$ to $\modu(k)$ which we shall continue to denote by
$R$; thus given $\phi$ a map in $\add(Q)$, its image under the functor
induced by $R$ is $R(\phi)$.
Let $\alpha$ be some dimension vector, and $\phi$ a map in $\add(Q)$
\begin{equation*}
\phi: \bigoplus_{v\in V} O(v)^{a(v)} 
\rightarrow \bigoplus_{v\in V} O(v)^{b(v)}
\end{equation*}
then for any representation $R$ of dimension vector $\alpha$, 
$R(\phi)$ is a ${\sum_{v\in V} a(v)\alpha(v)}$ by
${\sum_{v\in V}b(v)\alpha(v)}$ matrix.
If ${\sum_{v\in V}a(v)\alpha(v) = \sum_{v\in V}b(v)\alpha(v)}$, 
we define a semi--invariant 
polynomial function $P_{\phi,\alpha}$ on $R(Q,\alpha)$, by
\begin{equation*}
P_{\phi,\alpha}(p) = \det R_{p}(\phi).
\end{equation*}
We shall refer to these semi--invariants as determinantal
semi--invariants in future.  We will show  that the determinantal
semi--invariants   span all semi--invariants (Theorem
\ref{t2} below). To prove Thereom \ref{t2} we will use the classical
symbolic method which was also used by Procesi to show that the
invariants of matrices under conjugation are generated by traces
\cite{Procesi}. Procesi's result is generalized in \cite{LBP1} where
it is shown that invariants of quiver-representations are generated by
traces of oriented cycles. Subsequently Donkin showed that suitable
analogues of these results were valid in characteristic $p$ \cite{donkin,donkin1}.

It was shown in \cite{schofield} that the
determinantal semi-invariants can be defined in terms of the
representation theory of $Q$. This is reviewed in Section
\S\ref{representationtheory}. 

One corollary worth stating of this representation theoretic
interpretation is the following: define a point $p$ of
$R(Q,\beta)$ to be semistable if some non--constant semi--invariant
polynomial function does not vanish at $p$.  Then we have:
\begin{corollary} 
\label{interestingcorollary}
In characteristic zero
  the point $p$ of $R(Q,\beta)$ is semistable if and only there is some
  non--trivial (possibly infinite dimensional) representation $T$ of
  $Q$ such that
  $\Hom(T,R_p)=\Ext(T,R_p)=0$.
\end{corollary} 
For the proof we refer to the end of Section \ref{representationtheory}.

\medskip
In the sequel we will frequently change the quiver
 $Q=(V,A)$ to another quiver  $Q'=(V',A')$ which is connected to $Q$
 through an additive funcor.
\[
s:\add(Q')\rightarrow \add(Q)
\]
Through functoriality $s$ will act on various objects associated to
$Q$ and $Q'$. We will list these derived actions below. In order to avoid
having to introduce a multitude of adhoc notations we denote each of the
derived actions also by $s$.

To start there is an associated functor
\[
s:\Rep(Q)\rightarrow \Rep(Q'):R\mapsto R\circ s
\]
If $R$ has dimension vector $\alpha$ then $s(R)$ has dimension vector
$s(\alpha)\overset{\mathrm{def}}{=}\alpha\circ s$.  Put
$\alpha'=s(\alpha)$. We obtain that $s$ defines a map
\[
s:R(Q,\alpha)\rightarrow R(Q',\alpha')
\]
such that $s(R_p)=R_{s(p)}$. As usual there is
a corresponding $k$-algebra homomorphism
\[
s:S(Q',\alpha)\rightarrow S(Q,\alpha)
\]
given by $s(f)=f\circ s$. Writing out the definitions one obtains:
\begin{equation}
\label{relationwithdeterminantalinvariants}
s(P_{\phi',\alpha'})=P_{s(\phi'),\alpha}
\end{equation}
Finally $s$
defines a homomorphism
\[
s:\Gl(\alpha)\rightarrow \Gl(\alpha')
\]
which follows from functoriality by  considering $\Gl(\alpha)$ as the
automorphism group of $R_0(\alpha)$ in $\Rep(Q)$. One checks that for $g\in
\Gl(\alpha)$, $p\in R(Q,\alpha)$ one has $s(g\cdot p)=s(g)\cdot
s(p)$. It follows in particular that if $f$ is a semi--invariant in
$S(Q',\alpha')$ with character $\psi'$ then $s(f)$ is a semi-invariant
with character $s(\psi')\overset{\mathrm{def}}{=}\psi'\circ s$.

\section{Semi--invariant polynomial functions}
\label{s1}

Next, we discuss the ring $S(Q,\alpha)$.  $S(Q,\alpha)$ has two
gradings, one of which is finer than the other.  First of all,
$S(Q,\alpha)$ may be graded by $\mathbb{Z}^{A}$ in the natural way
since $R(Q,\alpha) = \times_{a\in A}{^{\alpha(ia)}k^{\alpha(ta)}}$.
We call this the $A$--grading.  On the other hand, $\Gl(\alpha)$ acts
on $R(Q,\alpha)$ and hence $S(Q,\alpha)$ and so $\times_{v\in V}
k^{*}$ acts on $S(Q,\alpha)$; we may therefore decompose $S(Q,\alpha)$
as a direct sum of weight spaces for the action of $\times_{v\in V}
k^{*}$ which gives a grading by $\mathbb{Z}^{V}$.  We call this the
$V$--grading.  The semi--invariants $P_{\phi,\alpha}$ are homogeneous
with respect to the second grading though not the first.  The first
grading is induced by the natural action of $\times_{a\in A} k^{*}$ on
$\add(Q)$. $\times_{a\in A} k^{*}$ acts on $\add(Q)$ by
$(\dots,\lambda_{a},\dots)(a) = \lambda_{a}a$ and according to
\eqref{relationwithdeterminantalinvariants} $g P_{\phi,\alpha} = P_{g \phi, \alpha}$ for $g \in \times_{a\in
  A} k^{*}$.

A standard Van der Monde determinant argument implies that the
vector subspace spanned by determinantal semi--invariants in
$S(Q,\alpha)$ is also the space spanned by the homogeneous components
of the determinantal semi--invariants with respect to the
$\mathbb{Z}^{A}$--grading. Thus it is enough to find the latter subspace.
Given a character $\chi$ of $\times_{a\in A} k^{*}$,
we define $P_{\phi,\alpha,\chi}$ to be the $\chi$--component of
$P_{\phi,\alpha}$.  

We begin by describing the semi--invariants which
are homogeneous with respect to the $A$--grading and linear in each
component of $R(Q,\alpha) = \times_{a\in
A}{^{\alpha(ia)}k^{\alpha(ta)}}$.  We call these the homogeneous
multilinear semi--invariants. We need temporarily another kind of
semi--invariant. A path $l$ in the quiver is an oriented cycle if $il =
tl$. Associated to such an oriented cycle is an invariant $\Tr_{l}$ for the action
of $\Gl(\alpha)$ on $R(Q,\alpha)$ defined by $\Tr_{l}(p) = \Tr(R_{p}(l))$
where $\Tr$ is the trace function.
\begin{theorem}
\label{t1}
The homogeneous multilinear semi--invariants of $R(Q,\alpha)$ are
spanned by semi--invariants of the form $P_{\phi,\alpha,\chi}
\prod_{i} \Tr_{l_{i}}$ where $l_{i}$ are oriented cycles in the quiver.
\end{theorem}

\begin{proof}
The semi--invariants are invariants for 
$\Sl(\alpha) = \times_{v \in V} \Sl_{\alpha(v)}(k)$ 
and conversely, $\Sl(\alpha)$--invariant 
polynomials that are homogeneous with respect to the $V$--grading are 
also semi--invariant for $\Gl(\alpha)$.
We may therefore use Weyl's 
description of the homogeneous multilinear invariants for $\Sl_{n}(k)$ 
and hence for $\Sl(\alpha)$.

Given a homogeneous 
multilinear $\Sl(\alpha)$--invariant
$f: \times_{a\in A}{^{\alpha(ia)}k^{\alpha(ta)}} \rightarrow k,$
it factors as
\begin{equation*}
\times_{a\in A}{^{\alpha(ia)}k^{\alpha(ta)}} \rightarrow 
\otimes_{a\in A} {^{\alpha(ia)}k^{\alpha(ta)}}
\xrightarrow{\tilde{f}} k
\end{equation*}
for a suitable linear map $\tilde{f}$.  Since 
\begin{equation*}
{^{\alpha(ia)}k^{\alpha(ta)}} \cong 
{^{\alpha(ia)}k}\otimes k^{\alpha(ta)}
\end{equation*}
as $\Sl(\alpha)$--representation, we may write 
\begin{equation*}
\tilde{f}: \otimes_{a\in A}{^{\alpha(ia)}k} \otimes k^{\alpha(ta)}
\rightarrow k.
\end{equation*}
Moreover, 
$\otimes_{a\in A}{^{\alpha(ia)}k} \otimes k^{\alpha(ta)}$ as 
$\Sl(\alpha)$--representation is a tensor product of covariant 
and contravariant vectors for $\Sl(\alpha)$.
Thus we may re--write
\begin{equation*}
\bigotimes_{a\in A} {^{\alpha(ia)}k} \otimes k^{\alpha(ta)} = 
\bigotimes_{v\in V} \left(\bigotimes_{a,ia=v} {^{\alpha(v)}k} 
\otimes \bigotimes_{a,ta=v} k^{\alpha(v)}\right)
\end{equation*}
and $\tilde{f} = \prod_{v\in V} \tilde{f}_{v}$ for 
$\Sl_{\alpha(v)}(k)$--invariant linear maps:
\begin{equation*}
\tilde{f}_{v}: \bigotimes_{a,ia=v} {^{\alpha(v)}k} 
\otimes \bigotimes_{a,ta=v} k^{\alpha(v)} \rightarrow k.
\end{equation*}
Roughly speaking Weyl \cite{weyl} showed that there are $3$ basic linear
semi--invariant functions on tensor products of covariant and
contravariant vectors for $\Sl_{m}(k)$. Firstly, there is the linear
map from ${^{m}k}\otimes k^{m}$ to $k$ given by $f(x\otimes y) = yx$;
this is just the trace function on ${^{m}k}\otimes k^{m} \cong
M_{m}(k)$. Secondly, there is the linear map from
$\otimes_{i=1}^{m}{^{m}k}$ to $k$ determined by $f(x_{1}\otimes \dots
\otimes x_{m}) = \det\left(x_1|\dots|x_m\right)$. The third case is
similar to the second; there is a linear map from
$\otimes_{i=1}^{m}k^{m}$ to $k$ again given by the determinant.
A spanning set for the linear semi--invariant functions on a general tensor
product of covariant and contravariant vectors is constructed from
these next.
A spanning set for the $\Sl_{m}(k)$--invariant 
linear maps from $\otimes_{B} {^{
m}k} \otimes \otimes_{C} k^{m}$ to $k$ is obtained in the 
following way.
We take three disjoint indexing sets $I,J,K$: we have surjective 
functions $\mu: B \rightarrow I \cup K, \: 
\nu: C \rightarrow J \cup K$ 
such that $\mu^{-1}(k)$ and 
$\nu^{-1}(k)$ have one element each for $k \in K$, and 
$\mu^{-1}(i)$ and $\nu^{-1}(j)$ have $m$ elements each for 
$i \in I, \: j \in J$.
We label this data by 
$\gamma = (\mu, \nu, I, J, K)$.  
To $\gamma$, we associate an $\Sl_{m}(k)$--invariant linear map
\begin{equation*}
f_{\gamma}\left(\bigotimes_{b \in B} x_{b} \otimes\bigotimes_{C} y_{c}\right) 
= \prod_{k\in K} y_{\nu^{-1}(k)}x_{\mu^{-1}(k)} 
\prod_{i\in I} \det\left(x_{b_{1}}|\dots|x_{b_{m}}\right) 
\prod_{j\in J} \det 
\begin{pmatrix}
y_{c_{1}}\\
\vdots\\
y_{c_{m}}
\end{pmatrix}
\end{equation*}
where $\{b_{1},\dots,b_{m}\} = \mu^{-1}(i), \; 
\{c_{1},\dots,c_{m}\} = \nu^{-1}(j)$.
Note that $f_{\gamma}$ is determined only up to sign since we 
have not specified an ordering of $\mu^{-1}(i)$ or $\nu^{-1}(j)$.

A spanning set 
for $\Sl(\alpha)$--invariant linear maps from
\begin{equation*}
\bigotimes_{v\in V} \left(\bigotimes_{a,ia = v}{^{\alpha(v)}k} 
\otimes \bigotimes_{a,ta=v} k^{\alpha(v)}\right)
\end{equation*}
is therefore determined by giving quintuples 
$(\mu,\nu,I,J,K) = \Gamma$ where
\begin{align*}
I &= \bigcup_{v\in V}^{\bullet} I_{v}\\ 
J &= \bigcup_{v\in V}^{\bullet} J_{v}\\
K &= \bigcup_{v\in V}^{\bullet} K_{v}
\end{align*}
and surjective maps 
\begin{align*}
\mu &: A \rightarrow I \cup K,\\
\nu &: A \rightarrow J \cup K
\end{align*}
where 
\begin{align*}
\mu(a) &\in I_{ia} \cup K_{ia},\\
\nu(a) &\in J_{ta} \cup K_{ta},
\end{align*}
$\mu^{-1}(k)$ and $\nu^{-1}(k)$ have one element each, $\mu^{-1}(i)$
and $\nu^{-1}(j)$ have $\alpha(v)$ elements each for 
$i \in I_{v}\text{ and } j\in J_{v}$. 
Then $\Gamma$ determines data $\gamma_{v}$ for each $v \in
V$ and we define 
\begin{equation*}
f_{\Gamma} = \prod_{v\in V} f_{\gamma_{v}}.
\end{equation*}

We show that these specific semi--invariants lie in the linear span of
the homogeneous components of  determinantal semi--invariants.
First, we treat the case where $K$ is empty.  Let $n = |A|$.  
We have two expressions for $n$:
\begin{equation*}
n = \sum_{v\in V} |I_{v}| \alpha(v) = \sum_{v\in V} |J_{v}| \alpha(v).
\end{equation*}
To each arrow $a$, we have 
a pair $(\mu(a), \nu(a))$ associated.  To this data, we associate 
a map in $\add(Q)$ in the following way.  We consider a map
\begin{equation*}
\Phi_{\Gamma}: \bigoplus_{v\in V} O(v)^{I_{v}} \rightarrow
\bigoplus_{v\in V} O(v)^{J_{v}}
\end{equation*}
whose $(i,j)$--entry is 
\begin{equation*}
\sum_{a,(\mu(a),\nu(a))=(i,j)}a.
\end{equation*}
Given $p \in R(Q,\alpha), \; R_{p}(\Phi_{\Gamma})$ is an 
$n$ by $n$ matrix which we may regard as a partitioned
matrix where the rows are indexed by $I$ and the columns by $J$, 
there are $\alpha(v)$ rows having index $i \in I_{v}$, 
$\alpha(v)$ columns having index $j \in J_{v}$ and the block having 
index $(i,j)$ is 
\begin{equation*}
\sum_{a:(\mu(a),\nu(a))=(i,j)} R_{p}(a).
\end{equation*}
We claim that $f_{\Gamma} = P_{\Phi_{\Gamma},\alpha,\chi}$ 
up to sign where $\chi((\lambda_a)_a)=\prod_a \lambda_a$. To prove
this it will be convenient 
to define a new quiver $Q'=(V',A')$ whose vertices are given by $I\cup J$
and whose edges are the same as those of $Q$. The initial and terminal
vertex of $a\in A'=A$ are given by $(\mu(a),\nu(a))$. On $Q'$ we define
data $\Gamma'$ which is defined by the same quintuple
$(\mu,\nu,I,J,K)$ as $\Gamma$, but with has different decompositions
$I=\bigcup_{i\in V'} I_{i}$, $J=\bigcup_{j \in V'} J_{j}$. In fact 
\[
I_{i}=\begin{cases}
\{i\}&\text{if $i\in I$}\\
\emptyset&\text{otherwise}
\end{cases}
\qquad
\text{and}
\qquad
J_{j}=\begin{cases}
\{j\}&\text{if $j\in J$}\\
\emptyset&\text{otherwise}
\end{cases}
\]
We define a functor $s:\add(Q')\rightarrow \add(Q)$ by $s(a)=a$,
$s(O(i))=O(v)$, $s(O(j))=O(w)$ where $a\in A'=A$, $i\in I_v$, $j\in
J_w$.

Since $Q$ and $Q'$ have the same edges, the action of $\times_{a\in
  A}k^\ast$ on $\add(Q)$ lifts canonically to an action of
$\times_{a\in A}k^\ast$ on $\add(Q')$. Put $\alpha'=s(\alpha)$. Using
\eqref{relationwithdeterminantalinvariants} we then find
$s(P_{\Phi_{\Gamma'},\alpha'})=P_{s(\Phi_{\Gamma'}),\alpha}
  =P_{\Phi_\Gamma,\alpha}$ and similarly
  $s(P_{\Phi_{\Gamma'},\alpha',\chi})=P_{\Phi_\Gamma,\alpha,\chi}$.
  Finally one also verifies that $s(f_{\Gamma'})=f_\Gamma$.

Hence to prove that $f_\Gamma=\pm P_{\Phi_\Gamma,\alpha,\chi}$ we may
  replace the triple $(Q,\alpha,\Gamma)$ by $(Q',\alpha',\Gamma')$. We
  do this now.

In order to prove that $f_\Gamma=\pm P_{\Phi_\Gamma,\alpha,\chi}$, we need only 
check that the two functions agree on the image of 
$W=\times_{a\in A}{^{\alpha(ia)}k} \times k^{\alpha(ta)}$ in 
$\otimes_{a\in A} {^{\alpha(ia)}}k \otimes k^{\alpha(ta)}$.

Let $\psi:\Gl(\alpha)\rightarrow k^\ast$ be the
character given by 
\[
\psi((A_v)_{v\in V})=\prod_{i\in I}\det A_i \cdot
\prod_{j\in J} (\det A_j)^{-1}
\]
Then one checks that both $f_\Gamma$ and $P_{\Phi_\Gamma,\alpha}$ are
semi-invariants on $W$ with character $\psi$. Now we claim that on $W$
we have $f_\Gamma=P_{\Phi_\Gamma,\alpha}$, up to sign. To prove this
we use the fact that $\Gl(\alpha)$ has an open orbit on $W$.

For  vertices $i \in I$ and  $j \in J$, we let 
$\{a_{i,1},\dots,a_{i,\alpha(i)}\}$ and
$\{a_{j,1},\dots,\alpha_{j,\alpha(j)}\}$ be  the sets of
arrows incident with $i$ and $j$ respectively.
So
\begin{equation*}
\times_{a\in A}{^{\alpha(ia)}k} \times k^{\alpha(ta)} = 
\times_{i\in I} \left(\times_{l=1}^{\alpha(i)}{^{\alpha(i)}k}\right) \times 
\times_{j\in J} \left(\times_{m=1}^{\alpha(j)} k^{\alpha(j)}\right).
\end{equation*}
We take the point $p$ whose $(i,l)$--th entry is the $l$--th standard
column vector in $^{\alpha(i)}k$; that is its $n$th entry is $\delta_{nl}$ 
and whose $(j,m)$--th entry is the $m$--th standard row vector in 
$k^{\alpha(j)}$. The $\Gl(\alpha)$-orbit of this point is open in $W$. Hence to
show that $f_\Gamma=\pm P_{\Phi_\Gamma,\alpha}$ on $W$, it suffices to do this
in the point $p$.

 Now $f_{\Gamma}(p)=\pm 1$. Furthermore we have
\[
R_{p}(\Phi_{\gamma})_{(i,l),(j,m)}=
\begin{cases}
1&\text{if $a_{il}=a_{jm}$}\\
0&\text{otherwise}
\end{cases}
\]
In particular 
 $R_p(\Phi_{\Gamma})$ is a 
permutation matrix and thus $P_{\Phi_\Gamma,\alpha}(p)=\pm 1$.
Hence indeed
$f_{\Gamma} = \pm P_{\Phi_{\Gamma},\alpha}$ (on $W$). Note that the
 non-zero entries of $R_p(\Phi_{\Gamma})$ are naturally indexed by
 $A$.

 To prove that $f_\Gamma=\pm P_{\Phi_\Gamma,\alpha,\chi}$ on $W$ it is now
 sufficient to prove that
 $P_{\Phi_\Gamma,\alpha}=P_{\Phi_\Gamma,\alpha,\chi}$ on $W$.
To this end we lift the $\times_{a\in A} k^\ast$  action on $R(Q,\alpha)$ to
 $W$ by defining $(\lambda_a)_a\cdot (x_{ia},y_{ta})_a=
 (\lambda_a x_{ia},y_{ta})_a$ (this is just some convenient
 choice). 

Now we have to show  that $P_{\Phi_\Gamma,\alpha}$ is itself
   homogeneous with character $\chi$ when restricted to $W$. Since the
   action of $\times_{a\in A} k^\ast$ commutes with the $\Gl(\alpha)$-action it
   suffices to do this in the point $p$. Now if we put $q=(\lambda_a)_a\cdot
   p$  then $R_q(\Phi_\Gamma)$ is obtained from $R_p(\Phi_\Gamma)$ by
    multiplying by $\lambda_a$ for all $a\in A$ the non-zero entry in
   $R_p(\Phi_\Gamma)$ indexed by $a$. Therefore $P_{\Phi_\Gamma,
     \alpha}(p)$ is multiplied by $\prod_a\lambda_a$. This proves what
   we want.

\medskip

It remains to deal with the case where $K$ is non--empty.  
Roughly speaking one of two things happens here; if we have two
distinct arrows $a=\nu^{-1}(k)$ and $b=\mu^{-1}(k)$ for some $k\in K$
then this element of $K$ corresponds to replacing $a$ and $b$ by their
composition; if on the other hand $a=\nu^{-1}(k)=\mu^{-1}(k)$ then
$ia=ta$ and this element of $K$ corresponds to taking the trace of $a$.

We associate a quiver $Q(A,K)$ with vertex set $A$ and arrow set $K$;
given $k \in K$, $ik = \nu^{-1}(k), \; tk = \mu^{-1}(k)$. This quiver
has very little to do with the quiver $Q$; it is a temporary
notational convenience.  The connected components of $Q(A,K)$ are of
three types: either they are oriented cycles, open paths or isolated
points. The vertices of components of the first type are arrows of $Q$
that also form an oriented cycle, those of the second type are arrows
that compose to a  path in the $Q$ (which can in fact also be an
oriented cycle); the isolated points we shall
treat in the same way to the second type.

We label the oriented cycles $\{L_{l}\}$ and the open paths $\{M_{m}\}$.
To an oriented cycle $L$ 
in $Q(A,K)$, we associate the invariant $\Tr_{p_{L}}$ where 
$p_{L}$ is the path around the oriented cycle in $Q$. This is
independent of our choice of starting point on the loop.
To the open path $M$ in $Q(A,K)$, we associate the path $p_{M}$,
the corresponding path in $Q$.

We consider the adjusted quiver $Q_{K}$ with vertex 
set $V$ and arrow set $\{p_{M_{m}}\}$ where $i{p_{M}}$ and 
$t{p_{M}}$ are defined as usual. Define the functor
$s:\add(Q_K)\rightarrow \add(Q)$ as follows: on vertices $s$ is the
identity, and on edges $s(p_{M_m})=p_{M_m}$. 

If $p_{M_{m}} = a_{m,1} \dots a_{m,d}$, we define
\begin{align*}
\mu(p_{M_{m}}) &= \mu(a_{m,1})\\
\nu(p_{M_{m}}) &= \nu(a_{m,d}).
\end{align*}

Then
\begin{align*}
\mu &: \{p_{M_{m}}\} \rightarrow I\\
\nu &: \{p_{M_{m}}\} \rightarrow J
\end{align*}
are surjective functions giving data $\Gamma_{K}$ on $Q_{K}$.
One checks directly that
\begin{align*}
f_{\Gamma} = s(f_{\Gamma_{K}}) \prod_{L} \Tr_{p_{L}}
&=s(P_{\phi_{\Gamma_{K}},\alpha,\chi_{K}})\prod_{L} \Tr_{p_{L}}\\
&=P_{s(\phi_{\Gamma_K}),\alpha,\chi}\prod_{L} \Tr_{p_{L}}
\end{align*}
for $\chi_{K}$ and $\chi$ of weight $1$ on each arrow for the quivers
$Q_{K}$ and $Q$ respectively, which completes our proof.
\end{proof}
If $l$ is an oriented cycle in the quiver then $\Tr(R_{p}(l))$ 
lies in the linear span of $\det(I + \lambda R_{p}(l))$ 
for varying $\lambda \in k$ and 
$\det(I + \lambda R_{p}(l)) = P_{I + \lambda l,\alpha}(p)$
Also 
\begin{equation*}
P_{\phi,\alpha} P_{\mu,\alpha} = 
P_{\bigl(\begin{smallmatrix}
\phi & 0\\
0 & \mu
\end{smallmatrix}\bigr),\alpha}
\end{equation*}
so we may deduce the following corollary.
\begin{corollary}
\label{t3}
Homogeneous multilinear semi--invariants lie in
the linear span of determinantal semi--invariants.
\end{corollary}
Since we have dealt with the multilinear case, standard 
arguments in characteristic $0$ apply to the general case.
\begin{theorem}
\label{t2}
In characteristic $0$, the semi--invariant polynomial functions for
the action of $\Gl(\alpha)$ on $R(Q,\alpha)$ are spanned by
determinantal semi--invariants.
\end{theorem}
\begin{proof}
It is enough to consider semi--invariants 
homogeneous with respect to the $A$--grading of weight $\chi$ 
where $\chi((\lambda_a)_a) = \prod_a \lambda_a^{m_{a}}$.

 We may also assume that no component of the
$A$--grade on the semi--invariant is $0$ since we may always restrict to
a smaller quiver.
To $\chi$, we associate a new quiver $Q_{\chi}$ with vertex set 
$V$ and arrow set $A_{\chi}$ where 
$A_{\chi} = \bigcup_{a\in A} \{a_{1},\dots,a_{m_{a}}\}$ 
where $ia_{i} = ia, \; ta_{i} = ta$.

We have functors
\begin{gather*}
\sigma: \add(Q) \rightarrow \add(Q_{\chi})\\
\sigma(a) = \sum a_{i}\\
\pi: \add(Q_{\chi}) \rightarrow \add(Q)\\
\pi (a_{i}) = a.
\end{gather*}

Given a semi--invariant $f$ for $R(Q,\alpha)$, we define 
$\tilde{f}$ to be the $\chi'$--component of $\sigma(f)$ where 
$\chi'((\lambda_{a_i})_{a,i}) = \prod_{a,i} \lambda_{a_i}$.
So $\tilde{f}$ is homogeneous multilinear. Then one checks that
$\pi(\tilde{f)}=\prod_a m_{a}!\cdot f$.

However, by the previous corollary, $\tilde{f}$ lies in the linear
span of determinantal semi--invariants; therefore
\begin{equation*}
\tilde{f}=\sum_{i}\lambda_{i}P_{\phi_{i},\alpha} 
\end{equation*}
and so
\begin{align*}
\prod_{a\in A} m_{a}! \cdot f &=
\pi(\tilde{f})\\
&=\sum_{i}\lambda_{i}\pi(P_{\phi_{i},\alpha})\\
&=\sum_{i}\lambda_{i}P_{\pi(\phi_{i}),\alpha}\qquad 
\text{(by \eqref{relationwithdeterminantalinvariants})} 
\end{align*}
Hence, in characteristic 0, semi--invariants homogeneous with respect
to the $A$--grading and hence all semi--invariants lie in the linear
span of the determinantal semi--invariants as required.
\end{proof}

To finish this section we will give a 
a slightly more combinatorial
interpretation
of our proof of Theorem \ref{t2}. 

Let us say that a pair $(Q,\alpha)$ is \emph{standard} if one of
the following holds.
\begin{enumerate}
\item $Q$ consists of one vertex with one loop.
\item $Q$ is bipartite with vertex set $I\cup J$ with all initial
  vertices in $I$ and all terminal vertices in $J$. Furthermore if
  $v\in V=I\cup J$ then $v$ belongs to exactly $\alpha(v)$ edges.
\end{enumerate}
To a standard pair $(Q,\alpha)$ we associate a standard semi--invariant
 as follows.
If $Q$ is a loop then we take $\Tr(R_p(a))$ where $a$ is the loop. If
$Q$ is bipartite then we take $P_{\phi,\alpha}$ where $\phi$ is given
by
\[
\phi:\oplus_{i\in I} O(i)\rightarrow \oplus_{j\in J} O(j)
\]
such that the $(i,j)$ component of $\phi$ is given by the sum of
arrows from $i$ to $j$. 

If $(Q,\alpha)$ is arbitrary then we define the standard
semi--invariants of $S(Q,\alpha)$ as those semi--invariants which are
of the form $s(f)$ for some functor
\[
s:\add(Q')\rightarrow \add(Q)
\]
such that $(Q',s(\alpha))$ is standard and such that $f$ is the
corresponding standard semi--invariant.
 
The follows corollary can now easily be obtained from our proof of
Theorem \ref{t2}.
\begin{corollary}
In  characteristic zero the semi--invariants in $S(Q,\alpha)$ are generated by the standard
semi--invariants. 
\end{corollary}

\section{Interpretation in terms of representation theory}
\label{representationtheory}
This section is mainly a review of some results of \cite{schofield} which
provided the motivation for the current paper. We prove Corollary
\ref{interestingcorollary}. 

It is standard that the category $\Rep(Q)$ is equivalent to the
category of right modules over $kQ$, the \emph{path algebra} of
$Q$. We will identify a representation with its corresponding $kQ$-module.

With every vertex $v\in Q$ corresponds canonically an idempotent
$e_v$ in $kQ$ given by the empty path. We denote by $P_v$ the
 representations
of $Q$ associated to $e_v kQ$. Clearly this is a projective object in
$\Rep(Q)$. 

For any representation $R$ we have $\Hom(P_{v},R)\cong R e_v=R(v)$
and furthermore  $\Hom(P_v,P_w)=e_w kQ e_v$. Now $e_w kQ e_v$ is a
vector space spanned by the paths having initial vertex $w$ and
terminal vertex $v$. Hence in fact
\[
\Hom_{\Rep(Q)}(P_v,P_w)=\Hom_{\add(Q)}(O(w),O(v))
\]
In other words if we denote by $\proj(Q)$ the additive category
generated by the $(P_v)_{v\in V}$ then $\proj(Q)$ is equivalent to the
opposite category of $\add(Q)$. 

Now recall the following:
\begin{lemma}
\label{projectiveslemma}
$\proj(Q)$ is equivalent to the category of
finitely generated projective $kQ$-modules.  
\end{lemma}
Only in the case of quivers
with oriented cycles, there is something to prove here; one must show
that the finitely generated projective modules $kQ$ of the quiver are
all direct sums of $P_v$'s. This may be proved in a similar way as the
fact that over a free algebra all projective modules are free \cite{cohn1}.

\medskip

Given a map
\begin{equation*}
\gamma: \bigoplus_{v\in V} P_{v}^{b(v)} 
\rightarrow \bigoplus_{v\in V} P_{v}^{a(v)},
\end{equation*}
we denote by 
\begin{equation*}
\hat{\gamma}: \bigoplus_{v\in V} O(v)^{a(v)} 
\rightarrow \bigoplus_{v\in V} O(v)^{b(v)}
\end{equation*}
the corresponding map in $\add(Q)$; similarly, for $\mu$ a map in 
$\add(Q)$, $\hat{\mu}$ is the corresponding map in 
$\proj(Q)$.

In order to link the semi--invariant polynomial functions
$P_{\phi,\alpha}$ with the representation theory of $Q$, we recall
some facts and definitions about the category of representations of a
quiver. First of all, $\Ext^{n}$ vanishes for $n>1$. Given dimension
vectors $\alpha$ and $\beta$, we define the Euler inner product by
\begin{equation*}
\langle \alpha,\beta\rangle =
\sum_{v\in V}\alpha(v)\beta(v)-\sum_{a\in
A}\alpha(ia)\beta(ta).
\end{equation*} 
If $R$ and 
$S$ are representations of dimension vector $\alpha$ and $\beta$
respectively, then
\begin{equation}
\label{eulerformula}
\dim\Hom(R,S)-\dim\Ext(R,S)=\langle\alpha,\beta\rangle.
\end{equation}

Given representations $R$ and $S$, we say that $R$
is left perpendicular to $S$ and that $S$ is right perpendicular to
$R$ if and only if $\Hom(R,S)=0=\Ext(R,S)$; 
Given a representation $R$ we
define the right perpendicular category to $R$, $R^{\perp}$, to be the
full subcategory of representations that are right perpendicular to
$R$ and the left perpendicular category to $R$, ${^{\perp}R}$ is defined
to be the full subcategory of representations that are left
perpendicular to $R$. It is not hard to show  that $R^{\perp}$ and ${}^\perp R$
are exact hereditary abelian subcategories of $\Rep(Q)$. In \cite{schofield} the first author
even shows that if $Q$ has no oriented cycles and $R$ has an open orbit in $R(Q,\alpha)$
then $R^{\perp}$ is given by the representations of a
quiver without oriented cycles and with $|V|-s$ vertices where $s$ is
the number of non-isomorphic indecomposable summands of $R$. A similar
result holds for ${}^\perp R$.

If $R$ and $S$ are finite dimensional and $R$ is left perpendicular to
$S$ then it follows from \eqref{eulerformula} that $\langle\ddim R,\ddim
S\rangle = 0$. The converse problem is interesting: suppose that we
have a representation $R$ of dimension vector $\alpha$ and a
dimension vector $\beta$ such that $\langle\alpha,
\beta\rangle =0$, what are the conditions on a point $p\in R(Q,\beta)$
in order that $R_p$ lies in $R^{\perp}$?

This problem was discussed and solved by the first author in
\cite{schofield} in the case that $Q$ has no oriented cycles. Let us 
assume this for a moment. We start with a minimal projective
resolution of $R$:
\[
0\rightarrow \bigoplus_v P_v ^{b(v)}\xrightarrow{\theta} \bigoplus_v
P_v^{a(v)}
\rightarrow R\rightarrow 0.
\]
By the above discussion $\theta=\hat{\phi}$ for some map $\phi:\bigoplus_v
O(v)^{a(v)}\rightarrow \bigoplus_v O(v)^{b(v)}$ in
$\add(Q)$.  Applying $\Hom(-,R_p)$ yields a long exact sequence
\begin{equation}
\label{longexactsequence}
0\rightarrow \Hom(R,R_p)\rightarrow \Hom(\oplus_v
P_v^{a(v)},R)\xrightarrow{\Hom(\hat{\phi},R_p)}  \Hom(\oplus_v
P_v^{b(v)},R)\rightarrow
\Ext(R,R_p)\rightarrow 0
\end{equation}
and clearly $\Hom(\hat{\phi},R_p)=R_p(\phi)$.  The condition
$\langle\alpha,\beta\rangle=0$ translates into $\sum_{v\in A}
((a(v)-b(v))\beta_v=0$, so that $R_p(\phi)$ is in fact represented by
a square matrix.

In \cite{schofield} the first
author defined $P_{R,\beta}(p)=\det R_p(\phi)$. It is not hard to see
that $P_{R,\beta}$ is independent of the choice of $\theta$ and
furthermore that it is a polynomial on $R(Q,\beta)$.

Hence if we define
\[
V(R,\beta)=\{p\in R(Q,\beta)\mid R\perp R_p\}
\]
then it follows from \eqref{longexactsequence} that 
\[
V(R,\beta)=\{p\in R(Q,\beta)\mid P_{R,\beta}(p)\neq 0\}
\]
In  other words, $V(R,\beta)$ is either trivial or the complement of a
hypersurface.

If there is an exact sequence in $\Rep(Q)$
\[
0\rightarrow R_1\rightarrow R \rightarrow R_2\rightarrow 0
\]
with $\langle\ddim R_1,\beta\rangle=\langle \ddim R_2,\beta\rangle=0$
then clearly $P_{R,\beta}=P_{R_1,\beta}P_{R_2,\beta}$. This provides some
motivation for the main result of \cite{schofield} which we state below.
\begin{theorem} Assume that $Q$ has no oriented cycles and let $S$ be
  a representation with dimension vector $\beta$ which has an open
  orbit in $R(Q,\beta)$. Let $R_1,\ldots,R_t$ be the simple objects of
  ${}^\perp S$. Then the ring of semi--invariants for $R(Q,\beta)$ is
  a polynomial ring in the generators $P_{R_i,\beta}$.
\end{theorem}
Now let us go back to the general case. Thus we allow that $Q$ has
oriented cycles. In this case $kQ$ may be infinite dimensional, so that
it is not clear how to define the minimal resolution of a
representation $R$. Therefore we will take the map $\phi:\bigoplus_v
O(v)^{a(v)}\rightarrow \bigoplus_v O(v)^{b(v)}$ as our fundamental object and we put
$P_{\phi,\beta}(p)=\det R_p(\phi)$ (provided $\sum_{v\in A}
(a(v)-b(v))\beta_v=0$). 

To make the link with the discussion above we need that $\hat{\phi}$
is injective. What happens if $\hat{\phi}$ is not injective? Then
there is a non-trival kernel
\[
0\rightarrow P\rightarrow \bigoplus_v P_v ^{b(v)}\xrightarrow{\hat{\phi}} \bigoplus_v P_v^{a(v)}
\]
By the fact that $kQ$ is hereditary and lemma \ref{projectiveslemma}
$P\cong\oplus_v P_v^{c(v)}$. Furthermore, again because $kQ$ is
hereditary $P$ is a direct summand of $\oplus_v P_v ^{b(v)}$. Now
using the equivalence with $\add(Q)$ we find that  $\hat{\phi}$ is
not injective if and only if $\oplus_v O(v) ^{b(v)}$ has a direct
summand in $\add(Q)$ which is in the kernel of
$\phi$.  Let $\phi'$ be the restriction of $\phi$ to the complementary
summand. Then we find
\[
P_{\phi,\beta}=
\begin{cases}
P_{\phi',\beta}&\text{if $\forall: v\in\Supp\beta:c(v)=0$}\\
0&\text{otherwise}
\end{cases}
\]

So for the purposes of semi-invariants we may assume that $\hat{\phi}$
is injective, which is what we will do below. If $Q$ has no oriented
cycles then it follows that $P_{\phi,\beta}=P_{\cok
  \hat{\phi},\beta}$. In general we 
 find (using
\eqref{longexactsequence}):
\begin{equation}
\label{property}
P_{\phi,\beta}(p) \neq 0 \Leftrightarrow 
\det R_{p}(\phi) \neq 0 \Leftrightarrow 
R_p \in {\cok\hat{\phi}}^{\perp} 
\end{equation}
Note however that $\cok \hat{\phi}$ may be infinite dimensional.
\begin{proof}[Proof of corollary \ref{interestingcorollary}]
Since the determinantal semi--invariant polynomial functions span all
semi--invariant polynomial functions, there must exist some $\phi\in
\add(Q)$ with the properties $P_{\phi,\beta}(p)\neq 0$,
$P_{\phi,\beta}$ is not constant and
 $\hat{\phi}$ is injective. Then
$\cok\hat{\phi}\in{^{\perp}R_{p}}$.  If $\cok\hat{\phi}=0$ then by
\eqref{property} $P_{\phi,\beta}$ is nowhere vanishing, and hence
constant by the Nullstellensatz. This is a contradiction, whence we
may take $T=\cok\hat{\phi}$.
\end{proof}

\begin{thebibliography}{1}

\bibitem{cohn1}
P.~M. Cohn, {\em Free rings and their relations.}, London Mathematical Society
  Monographs, vol.~19, Academic Press, Inc., 1985.

\bibitem{donkin}
S.~Donkin, {\em Invariants of several matrices}, Invent. Math. {\bf 110}
  (1992), no.~2, 389--401.

\bibitem{donkin1}
\bysame, {\em Polynomial invariants of representations of quivers}, Comment.
  Math. Helv. {\bf 69} (1994), no.~1, 137--141.

\bibitem{LBP1}
L.~Le~Bruyn and C.~Procesi, {\em Semisimple representations of quivers}, Trans.
  Amer. Math. Soc. {\bf 317} (1990), no.~2, 585--598.

\bibitem{Procesi}
C.~Procesi, {\em Invariant theory of {$n\times n$}-matrices}, Adv. in Math.
  {\bf 19} (1976), 306--381.

\bibitem{schofield}
A.~Schofield, {\em Semi-invariants of quivers}, J. London Math. Soc. (2) {\bf
  43} (1991), no.~3, 385--395.

\bibitem{weyl}
H.~Weyl, {\em The classical groups}, Princeton University Press, 1946.

\end{thebibliography}
\ifx\undefined\bysame
\newcommand{\bysame}{\leavevmode\hbox to3em{\hrulefill}\,}
\fi

\end{document}